\documentclass[reqno,oneside,10pt]{amsproc}

\usepackage{amsfonts,amsmath,amsthm,amscd,amssymb,latexsym}

\usepackage{graphicx,enumerate,xy,xspace}

\newcommand{\recuo}{\hspace{-0.05 cm}}

\newcommand{\prho}{\overline{\rho}} 
\newcommand{\can}{\beta} 
\newcommand{\tr}{\hat{t}} 
\newcommand{\um }{1_A} 

\newcommand{\hits}{\rhd}

\newcommand{\seta}{\rhd}
\newcommand{\rseta}{\lhd}
\newcommand{\hitby}{\leftharpoonup} 
\newcommand{\Hom}{\operatorname{Hom}}

\newcommand{\vai}{\rightarrow}

\newtheorem{lemma}{Lemma} 
\newtheorem{prop}{Proposition} 
 
\newtheorem{teo}{Theorem}
\newtheorem{thm}{Theorem}
\newtheorem{defi}{Definition}

 \title[Partial Hopf actions, partial invariants and a Morita context]{Partial Hopf actions, partial invariants and a Morita context}

 \author[M.M.S. Alves]{Marcelo \ Muniz \ S. \ Alves}
 \address{Departamento de Matem\'atica, Universidade Federal do Paran\'a, Brazil}
 \email{marcelomsa@ufpr.br}
 \author[E. Batista]{Eliezer Batista}
 \address{Departamento de Matem\'atica, Universidade Federal de Santa Catarina, Brazil}
 \email{ebatista@mtm.ufsc.br}


\begin{document}

\begin{abstract}

Partial actions of Hopf algebras can be considered as a generalization of partial actions of groups on algebras. Among important properties of partial Hopf actions, it is possible to assure the existence of enveloping actions \cite{paper}. This allows to extend several results from the theory of partial group actions to the Hopf algebraic setting. In this article, we explore some properties of the fixed point subalgebra with relations to a partial action of a Hopf algebra. We also construct, for partial actions of finite dimensional Hopf algebras a Morita context relating the fixed point subalgebra and the partial smash product. This is a generalization of a well known result in the theory of Hopf algebras \cite{susan} for the case of partial actions. Finally, we study Hopf-Galois extensions and reobtain some classical results in the partial case.
\end{abstract}
%
%
%
\maketitle

\section*{Introduction}

Partial group actions were first defined by R. Exel in the context of operator algebras and
they turned out to be a powerful tool in the study of $C^*$-algebras generated
by partial isometries on a Hilbert space \cite{ruy}. The developments
originated by the definition of partial group actions, soon became an independent topic of interest in ring theory \cite{dok}. Now, the results are formulated in a purely 
algebraic way, independent of the  $C^*$ algebraic techniques which originated them. 

A partial action $\alpha$ of a group $G$ on a (possibly non-unital) $k$-algebra $A$ is a pair of families of subsets of $A$ and maps indexed by $G$, \linebreak $\alpha = (\{\alpha_g\}_{g \in G}, \{D_g\}_{g \in G})$, where each $D_g$ is an ideal of $A$ and each $\alpha_g$ is an algebra isomorphism $\alpha: D_{g^{-1}} \rightarrow D_g$ satisfying the following conditions:
\begin{enumerate}
\item[(i)] $D_e = A$ and $\alpha_e = I_A$;
\item[(ii)] $\alpha_g(D_{g^{-1}} \cap D_h) = D_g \cap D_{gh}$ for every $g,h \in G$;
\item[(iii)] $\alpha_g(\alpha_h(x)) = \alpha_{gh}(x)$ for every 
$x \in D_{g^{-1}} \cap D_{(gh)^{-1}}.$ 
\end{enumerate}

A first example of partial action is the following: If $G$ acts on a 
algebra $B$ by automorphisms and $A$ is an ideal of $B$, then we
have a partial action $\alpha$ on $A$ in the following manner: letting
$\beta_g$ stand for the automorphism corresponding to $g$, take 
$D_g = A \cap \beta_g(A)$, and define $\alpha_g: D_{g^{-1}} \vai D_g$ 
as the restriction of the automorphism $\beta_g $ to $D_g$. 

Partial Hopf actions were motivated by an attempt to
generalize the notion of partial Galois extensions of commutative rings \cite{paques} to a broader context. The definition of partial Hopf actions and co-actions were introduced by S. Caenepeel and K. Janssen in \cite{caen06}, using the notions of partial entwining structures. In particular, partial actions of $G$ determine 
partial actions of the group algebra $kG$ in a natural way. In the
same article, the authors also introduced the concept of partial 
smash product, which in the case of the group algebra $kG$, turns 
out to be the crossed product by a partial action 
$A\rtimes_{\alpha} G$. Further developments in the theory of
partial Hopf actions were done by C. Lomp in \cite{lomp07}, where the author
pushed forward classical results of Hopf algebras concerning smash products, 
like the Blattner-Montgomery and Cohen-Montgomery theorems \cite{susan}.

In \cite{paper}, we proved the theorem of existence of an enveloping action for a partial Hopf action, that is, if $H$ is a Hopf algebra which acts partially on a unital algebra $A$, 
then there exists an $H$-module algebra $B$ such that $A$ is isomorphic 
to a right ideal of $B$, and the restriction of the action of $H$ to this ideal is 
equivalent to the partial action of $H$ on $A$. Basically, the same ideas for the proof of the existence of an enveloping action for a partial group action \cite{dok} are present in the Hopf algebraic case. In the same article, we also proved many results related to the enveloping action: the existence of a Morita context between the partial smash product $\underline{A\# H}$, where $H$ is a Hopf algebra which acts partially on the unital algebra $A$, and the smash product $B\# H$, where $B$ is an enveloping action of $A$. The conditions for the existence of an enveloping co-action associated to a partial co-action of a Hopf algebra $H$ on a unital algebra $A$. Finally, we introduced the notion of partial representation of a Hopf algebra and showed that, under certain conditions on the algebra 
$H$, the partial smash product $\underline{A\# H}$ carries a partial representation of
$H$. 

In this work, we push forward the results obtained in \cite{paper}. First, we define the invariant sub-algebra $A^{\underline{H}} \subseteq A$, where $H$ is a Hopf algebra acting partially on a unital algebra $A$, and explore some of its properties. In what follows, we prove the existence of a Morita context between the invariant sub-algebra $A^{\underline{H}}$ and the partial smash product $\underline{A\# H}$, generalizing a classical result in the theory of Hopf algebras 
\cite{susan}. We also study Hopf-Galois extensions when $H$ is finite-dimensional and reobtain, in the partial case, classical results such as: if the trace mapping is surjective, then the Morita context is strict if and only if the extension $A^{\underline{H}} \subset A$ is Hopf-Galois. 

It is worth mentioning that in the paper \cite{caen06}, the authors considered a partial coaction of a Hopf algebra $H$ on a unital algebra $A$ and established a Morita context relating the subalgebra of coinvariants $A^{\underline{CoH}}$ and the dual smash product $\underline{\#}(H,A)$; here we use the fact that $H$ is finite-dimensional to build the Morita context relating directly $A^{\underline{H}}$ and $\underline{A \# H}$, and to study Hopf-Galois extensions with a more elementary approach.

\section{Partial Hopf actions}

We recall that a left action of a Hopf algebra $H$ on a unital algebra $A$ is  
a linear mapping $\alpha: H \otimes A \vai A$, which we will denote by 
$\alpha(h \otimes a) = h \hits a$, such that

\begin{enumerate}
\item[(i)] $h \hits (ab) = \sum (h_{(1)} \hits a) (h_{(2)} \hits b),$
\item[(ii)] $1 \hits a =a$
\item[(iii)]  $h \hits (k \hits a) = hk \hits a$ 
\item[(iv)] $h \hits \um = \epsilon(h) \um. $
\end{enumerate}
We also say that $A$ is an $H$ {\it module algebra}. Note that (ii) and (iii) say that $A$ is a left $H$-module.

\begin{defi} A partial action of the Hopf algebra $H$ on the 
algebra $A$ is a linear mapping $\alpha: H \otimes A \vai A$, denoted 
here by $\alpha(h \otimes a) = h \cdot a$, such that  
\begin{enumerate}
\item[(i)] $h \cdot (ab) = \sum (h_{(1)} \cdot a) (h_{(2)} \cdot b),$
\item[(ii)] $1 \cdot a = a,$
\item[(iii)] $h \cdot (g \cdot a) = \sum (h_{(1)} \cdot \um) ((h_{(2)} g) \cdot a).$
\end{enumerate}
In this case, we call $A$ a partial $H$ module algebra.
\end{defi} 

Because of the good dual properties of Hopf algebras, one can also define the concept of partial (right) coaction of a Hopf algebra $H$ on a unital algebra $A$, which will be important later.

\begin{defi} A partial coaction of the Hopf algebra $H$ on the 
algebra $A$ is a linear mapping $\prho : A \rightarrow H \otimes A $, denoted by 
$\prho (a)=\sum a^{[0]} \otimes a^{[1]}$, such that
\begin{enumerate}
\item[(i)] $\prho (a.b)=\prho (a) . \prho (b)$, $\forall a,b\in A$,
\item[(ii)] $(I\otimes \epsilon )\prho (a)=a$, $\forall a\in A$,
\item[(iii)] $(\prho \otimes I) \prho (a) =(\prho(\um )\otimes 1_H )
((I\otimes \Delta)\prho(a))$.
\end{enumerate}
\end{defi}

It is easy to see that every action is also a partial action (and the same holds for coactions).

As a basic example, consider  a partial action $\alpha$ of a group $G$ on an
unital algebra $A$. Suppose that each $D_g$ is also a unital algebra, that is,
$D_g$ is of the form $D_g =A1_g$ then there is a partial action of the group
algebra $kG$ on $A$ defined on the elements of the basis by
\begin{equation}
\label{actiongroupalgebra}
g\cdot a =\alpha_g (a 1_{g^{-1}}),
\end{equation}
and extended linearly to all elements of $kG$. 

There is an important class of examples of
partial Hopf actions: these induced by total actions. Basically, the induced partial actions can be described by the following result.

\begin{prop} \label{induced}
Let $H$ be a Hopf algebra, $B$ a $H$-module algebra and let $A$ be a
right ideal of $B$ with unity $\um $. 
Then $H$ acts partially on $A$ by 
\[
h \cdot a = \um  (h \hits a) 
\]
\end{prop}

The proof of this proposition is a straightforward calculation and can be seen in \cite{paper}.

If the algebra $A$ is a bilateral ideal of $B$, then its unit $\um$ is a central idempotent of $B$ and the induced partial action has a symmetric formulation, satisfying the additional relation \cite{paper},\\
{\ }\\
(iv) $h \cdot (g \cdot a) = \sum ((h_{(1)} g) \cdot a)(h_{(2)} \cdot \um).$\\
{\ }\\
This property is satisfied by partial actions of group algebras $kG$.

As a nontrivial example of partial Hopf action, we can consider the restriction of the action of the dual group algebra $kG^\ast$ of a finite group $G$ on the group algebra $kG$. Let $\{ p_g \}_{g\in G}$ be the dual basis for $kG^\ast$, the (global) action of $kG^\ast$ on $kG$ is given by $p_g \triangleright h=\delta_{g,h} h$. Consider now a normal subgroup $N\trianglelefteq G$, $N\neq \{ 1\}$, such that $\mbox{char}(k) \nmid |N|$; let $e_N \in kG$ be the central idempotent
\[
e_N =\frac{1}{|N|}\sum_{n\in N} n, 
\]
and let  $A$ be the ideal $A=e_N kG$, which is also a unital algebra with $\um = e_N$ . It is possible to restrict the action of $kG^\ast$ on $kG$ to a partial action on $A$. Given $x \in G$ and $p_g \in kG^*$, note that
\[
\begin{array}{rcl}
p_g \triangleright (e_N x ) & = & \sum_{h \in G}(p_{gh^{-1}} \triangleright e_N)
(p_h \triangleright x) \\
& = &  \sum_{h \in G}(p_{gh^{-1}} \triangleright e_N)p_h (x)x \\
& = &  (p_{gx^{-1}} \triangleright e_N)x\\
& = & \frac{1}{|N|} \sum_{n \in N}p_{gx^{-1}} (n) nx
\end{array}
\]
which is equal to $(1/|N|)g$ if $gx^{-1} \in N$, and is zero otherwise. Hence, if 
$gx^{-1} \in N$, 
 \[
 p_g \cdot (e_N x) = e_N (p_g \triangleright e_Nx) =  (1/|N|)e_Ng = (1/|N|)e_Nx
 \]
and $p_g \cdot (e_N x) =0$ otherwise. Therefore, 
$p_g \cdot e_N \neq \epsilon(p_g)e_N$ when $g \in N$ and the action is really partial.

On the other hand, given a partial action of a Hopf algebra $H$ on a unital algebra $A$, we can construct its enveloping action. For this intent, we need some preliminary definitions.

\begin{defi} Let $A$ and $B$ be two partial $H$-module algebras. We 
will say that a morphism of algebras $\theta: A \vai B$ is a morphism of 
partial $H$-module algebras if $\theta(h \cdot a) = h \cdot \theta (a)$ 
for all $h \in H$ and all $a \in A$. If $\theta$ is an isomorphism,
we say that the partial actions are equivalent.
\end{defi}

\begin{defi} Let $B$ be an $H$-module algebra and let $A$ be a right ideal of $B$ with unity $\um$. We will say that the induced partial action on $A$ is admissible  if $B = H \hits A$.
\end{defi}

\begin{defi}  Let  $A$ be a partial $H$-module algebra. An enveloping action for $A$ is a pair $(B,\theta)$, where 
\begin{enumerate}
\item[(i)] $B$ is a  (not necessarily unital) $H$-module algebra, that is, the item (iv) of the definition of total action doesn't need to be satisfied.
\item[(ii)] The map $\theta: A \vai B$ is a monomorphism of algebras.
\item[(iii)] The sub-algebra $\theta(A)$ is a right ideal in $B$.
\item[(iv)] The partial action on $A$ is equivalent to the induced partial action on $\theta (A)$.
\item[(v)]  The induced partial action on $\theta(A)$ is admissible.
\end{enumerate}
\end{defi}

Then we have the following result \cite{paper}.

\begin{teo}\label{existence}
Let $A$ be a partial $H$-module algebra and let $\varphi: A \vai Hom_k(H,A)$ be the map given by $\varphi(a)(h) = h \cdot a$, and let $B = H \hits \varphi(A)$; then $(B, \varphi)$ is an enveloping action of $A$.
\end{teo}

This enveloping action is called standard enveloping action, it is minimal in the sense that for every other enveloping action $(B', \theta)$, there is a epimorphism of $H$ module algebras $\Phi :B' \rightarrow B$ such that $\Phi (\theta (A))=\varphi (A)$.

As an example of globalization of partial Hopf actions, let $H_4$ be the Sweedler $4$-dimensional Hopf algebra, with $\beta = \{1,g,x,xg \}$ 
as a basis over the field $k$, where $char(k) \neq 2$. Another basis for $H_4$ is given by the elements
\[
e_1 = (1+g)/2,e_2 = (1-g)/2, h_1 = xe_1, 
h_2 = xe_2
\]
where the $e_i$'s form a complete system of primitive orthogonal idempotents of $H_4$, and the ideal generated by the  $h_i$'s is the radical of $H_4$ (the structure constants for the product and coproduct with respect to this basis can be found in \cite{paper}, for instance). Taking the dual basis $\beta^\ast = \{1^\ast, g^\ast, x^\ast, (xg)^\ast\}$, we obtain an isomorphism of Hopf algebras $\psi: H_4^\ast \vai H_4$ given by $1^\ast \mapsto e_1$, $g^\ast \mapsto e_2$,  $x^\ast \mapsto h_1$ and $ (xg)^\ast \mapsto h_2$. In the reference \cite{caen06} the authors constructed a partial coaction $\prho: k \vai k \otimes H_4$, given by $1 \mapsto 1 \otimes e$, where $e=\frac{1}{2}+\frac{1}{2}g-\alpha xg$. This partial coaction induces a partial $H_4^\ast $ action on $k$ by $h \cdot 1 = h(e)$ for each $h \in H_4^\ast$. For the basis elements of $H_4^\ast$, considering the above isomorphism, we have 
\[
e_1 \cdot 1 = \frac{1}{2}, \quad  e_2 \cdot 1 = \frac{1}{2}, \quad h_1 \cdot 1 = 0, \quad h_2 \cdot 1 = - \alpha. 
\]
Globalizing, we have the map $\varphi: k \vai \Hom_k(H_4^\ast, k) \cong H_4$ defined by $\varphi(1)( h)=(h \cdot 1)$, for each $h \in H_4^\ast$. Then, identifying $H_4^\ast$ with $H_4$ by $\psi$, we have  
\[
e_1 \rhd \varphi(1) = \frac{1}{2},  \quad e_2 \rhd \varphi(1) = \frac{1}{2}g - \alpha xg, \quad h_1 \rhd \varphi(1) = 0, \quad h_2 \rhd \varphi(1) = -\alpha 1  
\]
and hence the minimal enveloping action is given by $(B, \varphi)$, where $B = \langle 1, e \rangle_k$ (which is isomorphic to $k \times k$ as an algebra). 

In what follows, unless stated otherwise, we will consider only the case that $A$ is a bilateral ideal of $B$, and then the partial action satisfies the symmetric property (iv) above.

\section{Partial invariant subalgebras}

In the theory of partial group actions, one can define the invariant subalgebra in the following manner \cite{paques}: If $\alpha$ is a partial action of a group $G$ on a unital algebra $A$, such that each ideal $D_g$ is unital, for every $g\in G$, then the invariant subalgebra is the set
\[
A^{\alpha} =\{ a\in A\, |\, \alpha_g (a1_{g^{-1}}) =a1_g, \; \forall g\in G \} .
\]
It is an easy calculation to verify that $S^{\alpha}$ is, indeed, a sub-algebra of $A$. Motivated by this definition, we can define the invariant subalgebra by a partial action of a Hopf algebra. 

\begin{defi}
Let $H$ be a Hopf algebra acting partially on a $k$-algebra $A$. We define the set of invariants of the partial action as
\[ 
A^{\underline{H}} = \{a \in A; h \cdot a = a (h \cdot \um)  \} .
\]
\end{defi}

It is easy to prove that $A^{\underline{H}}$ is a subalgebra of $A$. Indeed, take $a,b\in A^{\underline{H}}$ and $h\in H$, then
\begin {eqnarray}
h\cdot (ab) &=& \sum (h_{(1)}\cdot a)(h_{(2)}\cdot b)=\sum a(h_{(1)}\cdot \um )(h_{(2)}\cdot b)=\nonumber\\
&=& a(h\cdot b) =ab(h\cdot \um ) .\nonumber
\end{eqnarray}

In the case of partial group actions each idempotent $1_g$ is central. More generally, it can be shown that if $H$ is cocommutative then  each element $h \cdot \um$ is central in $A$. Throughout this paper, {\it we assume that} $h \cdot \um$ {\it lies in the center of }$A$ {\it for each }$h \in H$.

In the case of a partial coaction of a Hopf algebra $H$ acting on a $k$-algebra $A$, we have also the notion of subalgebra of coinvariants.

\begin{defi}
Let $H$ be a Hopf algebra coacting partially on a $k$-algebra $A$. We define the set of coinvariants of the partial coaction as
\[ 
A^{\underline{CoH}} = \{a \in A |\,  \prho(a) = a \prho(\um )  \}
\]
\end{defi}
 
 In the same manner, one can easily verify that $A^{\underline{CoH}}$ is a subalgebra of $A$.

Let us consider now an enveloping action $(B, \theta)$ of the partial action of $H$ on $A$, or enveloping action of $A$, for short. If $B^H$ is the invariant subalgebra of $B$ with relation to the total action $\triangleright$ of $H$ on $B$, that is
\[
B^H =\{ b\in B \, |\, h\triangleright b =\epsilon (h)b\; \forall h\in H \} ,
\]
it is easy to see that $\um B^H \subseteq A^{\underline{H}}$, obviously, considering $\um B^H =\theta (\um )B^H$ and $A^{\underline{H}}=\theta (A^{\underline{H}})$. That is because if $b \in B^H$ then 
\begin{eqnarray}
h \cdot \um b &=& \um (h \triangleright \um b) = \um (\sum (h_{(1)} \triangleright \um )
(h_{(2)} \triangleright b)) =\nonumber \\
&=& \um (\sum (h_{(1)} \triangleright \um )(\epsilon(h_{(2)}) b)) =  
\um (h \triangleright \um) b = (h \cdot \um ) (\um b) \nonumber  
\end{eqnarray} 
On the other hand, the conditions to be satisfied in order to fulfill the equality are more restrictive, as we shall see later.

From now on, we will consider only the case where the Hopf algebra $H$ is finite dimensional and the partial action is symmetric, i.e., where property (iv) holds. 
Since $H$ is finite dimensional, there exists a nonzero left integral $t\in H$. Define the partial trace map $\tr: A \vai A$ by $\tr(a) = t \cdot a$. In the case of total actions this trace reduces to the classical one (see reference \cite{susan}, definition 4.3.3). The following two results are quite analogous to the classical results in the theory of Hopf algebras \cite{CFM}.

\begin{lemma}\label{trace}
$\tr$ is a $A^{\underline{H}} - A^{\underline{H}}$ bimodule mapping from $A$ into $A^{\underline{H}}$.
\end{lemma}

\begin{proof}
First, one need to check that $\tr (A) \subseteq A^{\underline{H}}$. Since $t$ is a left integral, given $a \in A$ and $h \in H$, 
\[
h \cdot \tr (a) = h \cdot ( t \cdot a) = 
\sum (h_{(1)} \cdot \um) (h_{(2)} t \cdot a) = \]
\[ = \sum (h_{(1)} \cdot \um) (\epsilon(h_{(2)}) t \cdot a) = (h \cdot \um ) \tr (a)
\]
and therefore $\tr(A) \subset A^{\underline{H}} $. 

Now the bimodule morphism property. If $b \in A^{\underline{H}}$ and $a\in A$, then
\[
\tr(ab) = t \cdot ab = \sum (t_{(1)} \cdot a)(t_{(2)} \cdot b) = \sum (t_{(1)} \cdot a) (t_{(2)} \cdot \um )b = (t \cdot a)b
\]

and
\[
\tr(ba) =  t \cdot ba = \sum (t_{(1)} \cdot b)(t_{(2)} \cdot a) = \sum  b (t_{(1)} \cdot \um)(t_{(2)} \cdot a) =  b (t \cdot a) =  b\tr(a) .
\] 
Note that in the third equality in the verification above we used the symmetry property (iv).
\end{proof}

\begin{prop} If $H$ is semisimple and unimodular, and $\um B^H = A^{\underline{H}}$, then the partial trace mapping is surjective. 
\end{prop}

\begin{proof} In fact, if $a \in A^{\underline{H}} $ then $a = \um b$ where $b \in B^{\underline{H}}$; since $H$ is semisimple, we have a nonzero left integral $t \in H$ such that $\epsilon(t) = 1$, and hence $\um b = \um \epsilon(t) b = \um (t \rhd b)$ (because $b \in B^H$). Now $b = \sum h_i \rhd a_i$ for some $h_i \in H$ and $a_i \in A$, and therefore, since $t$ is also a right integral, 
\[
a = \um (t \rhd b) = \um (t \rhd (\sum h_i \rhd a_i)) = \]
\[ = \um (\sum \epsilon(h_i) t \rhd a_i) = t  \cdot (\sum \epsilon(h_i)a_i)  \in t \cdot A. 
\]
\end{proof}

\begin{prop}
If the partial trace mapping is surjective then $\um B^H = A^{\underline{H}}$.
\end{prop}

\begin{proof} Let $x\in A^{\underline{H}}$. Since the partial trace is surjective, there is an element $y\in A$ such that $x=t\cdot y$, then
\[
x=t\cdot y=\um (t\triangleright y)\in \um B^H
\]
\end{proof}

Now, some words about the partial smash product. Let $A$ be a partial $H$ module algebra, we can endow the tensor product $A\otimes H$ with an associative algebra structure by
\[
(a\otimes h)(b\otimes k)=\sum a(h_{(1)} \cdot b) \otimes h_{(2)} k .
\]
Define the partial smash product as the algebra
\[
\underline{A\# H} =(A\otimes H)(1_A \otimes 1_H ) ,
\]
in other words, the smash product is the subalgebra generated by elements of the form
\[
\underline{a\# h} =\sum a(h_{(1)} \cdot \um )\otimes h_{(2)} .
\]
One can easily verify that the product with the symbol $\#$ satisfies
\[
(\underline{a\# h})(\underline{b\# k}) =\sum \underline{a(h_{(1)} \cdot b )\# h_{(2)}k} .
\]
Note that the elements of $A$ can be embedded into $\underline{A\# H}$ by the map 
$a\mapsto \underline{a\# 1_H }=a\otimes 1_H$, this is an algebra map. On the other hand, the elements of $H$ can be written into  the smash product as $\underline{\um \# h} =\sum (h_{(1)} \cdot \um )\otimes h_{(2)}$, but these elements do not form a subalgebra of $\underline{A\# H}$. Nevertheless, just as in the global case, the partial action is now implemented internally, since
\[
\sum \underline{(\um \# h_{(1)})} \ \ \underline{(a \# 1)} \ \ \underline{(\um \# S(h_{(2)}))} = \underline{(h \cdot a) \# 1}.
\]

The partial smash product has a very interesting factorization property, which will be useful later.

\begin{prop} Let $H$ be a Hopf algebra with invertible antipode acting partially on the algebra $A$. Then $\underline{A \# H} = (\um \otimes H )(A \otimes 1)$. 
\end{prop}

\begin{proof} In fact, consider $\underline{a\# h}\in \underline{A\#H }$, then we have
\begin{eqnarray}
\underline{a\# h} &=& \sum a(h_{(1)} \cdot \um ) \otimes h_{(2)} =\sum \epsilon (h_{(1)}) a (h_{(2)} \cdot \um )\otimes h_{(3)} =\nonumber\\
&=& \sum (h_{(2)}S^{-1}(h_{(1)})\cdot a)(h_{(3)} \cdot \um )\otimes h_{(4)} = \nonumber\\
& = &  
\sum h_{(2)}\cdot (S^{-1}(h_{(1)})\cdot a)\otimes h_{(3)} =\nonumber\\
&=& \sum (\underline{\um\# h_{(2)}})(\underline{(S^{-1} (h_{(1)})\# 1_H}) .\nonumber
\end{eqnarray}
\end{proof}

\begin{prop}
If $\tr: A \vai A^{\underline{H}}$ is surjective, then there is a non-zero idempotent $e$ in $\underline{A \# H}$ such that $e (\underline{A \# H}) e = (A^{\underline{H}} \# 1_H ) e \equiv A^{\underline{H}}$ as $k$-algebras.
\end{prop}

\begin{proof}
Here, for sake of simplicity, if $h \in H$ and $a \in A$, we will write $a$ for 
$\underline{a \# 1_H}$ and $h$ for $\underline{\um \# h}$. With this notation 
$ah=\underline{a\# h}$ and $ha=\sum (h_{(1)}\cdot a)h_{(2)}$, for every $a\in A$ and $h\in H$.

First, $hat = (h \cdot a)t$ for any $h \in H$ and $a \in A$:
\begin{eqnarray} 
hat &=& (\underline{\um \# h})(\underline{a \# 1_H}) (\underline{\um \# t}) = 
\sum (\underline{\um \# h})(\underline{ a \# t}) =\nonumber\\
&=& \sum \underline{(h_{(1)} \cdot a) \# h_{(2)}t} = 
\sum \underline{(h_{(1)} \cdot a)\# \epsilon(h_{(2)})t} = \nonumber\\
&=& \underline{h\cdot a) \#  t} = (h \cdot a) t .\nonumber
\end{eqnarray}

Assuming the surjectivity of the trace, let $c \in A$ be such that $\tr(c) = \um$, and consider $e = tc \in \underline{A \# H}$. This element is an idempotent, since
\[
e^2 = tctc = (tct)c = (t \cdot c) tc = \underline{\um \# 1_H} ( tc ) = tc = e.
\]

Let us verify that $e (\underline{A \# H}) e = (A^{\underline{H}} \# 1_H ) e$:
\begin{eqnarray}
e (ah) e & = & tcahtc  =   tca h \um t c =\nonumber \\
& = &  t(ca (h \cdot \um)) t c =  t \cdot (ca (h \cdot \um)) tc = \nonumber \\ 
& = & t \cdot (ca(h \cdot \um ))e =   \tr (ca(h\cdot \um ))e , \nonumber
\end{eqnarray}
which lies in $A^{\underline{H}} e$. 

Conversely, if $a \in A^{\underline{H}}$, then
\begin{eqnarray}
t \cdot (ca) &=& \sum (t_{(1)} \cdot c)(t_{(2)}\cdot a)= \nonumber\\
&=& \sum (t_{(1)} \cdot c)(t_{(2)}\cdot \um )a= \nonumber\\
&=& (t \cdot c)a = a,\nonumber 
\end{eqnarray}
and 
\[
ae = (t\cdot (ca))tc = tcatc = eae.
\]
Hence $e A^{\underline{H}} e = A^{\underline{H}} e.$

Finally, $e A^{\underline{H}} e$ is isomorphic to $A^{\underline{H}}$,  since for $a,b\in A^{\underline{H}}$
\[
(ae)(be) = at(cb)tc = a(t \cdot cb) tc = a(t \cdot c)btc = abtc = abe. 
\]
(where we used $t \cdot c = \um$).
\end{proof}

\section{A Morita context}

In what follows we show that, just as in the case of global actions (see reference \cite{susan} paragraphs 4.4 and 4.5), there is a Morita context connecting the algebras $A^{\underline{H}}$ and $\underline{A \# H}$. From now on, we are assuming that $H$ is a Hopf algebra with invertible antipode. Let us begin by recalling the definition of a Morita context between two rings. 

\begin{defi} A Morita context is a six-tuple $(R,S,M,N,[\cdot , \cdot ], \langle \cdot , \cdot \rangle )$ where
\begin{enumerate}
\item $R$ and $S$ are rings,
\item $M$ is an $R-S$ bimodule, 
\item $N$ is an $S-R$ bimodule, 
\item $[ \cdot , \cdot ]:M\otimes_S N\rightarrow R$ is a bimodule morphism, 
\item $ \langle \cdot , \cdot \rangle :N\otimes_R M \rightarrow S$ is a bimodule morphism, 
\end{enumerate}
such that
\begin{equation}
\label{assoctau}
[ m, n] m'=m \langle n,  m' \rangle , \qquad \forall m, m' \in M, \quad
\forall n \in N , 
\end{equation}
and 
\begin{equation}
\label{assocsigma}
\langle n,  m \rangle  n' = n  [m, n' ], \qquad \forall m\in M, \quad
\forall n, n' \in N.
\end{equation}
\end{defi}

By a fundamental theorem due to Morita (see, for example \cite{jacobson} on
pages 167-170), if the morphisms $[,]$ and $\langle, \rangle$ are surjective, then the
categories ${}_R\! \mbox{Mod}$ and ${}_S\! \mbox{Mod}$ are equivalent. In this
case, we say that $R$ and $S$ are Morita equivalent.

For the Morita context between $A^{\underline{H}}$ and $\underline{A \# H}$, the bimodules $M$ and $N$ will both have $A$ as subjacent vector space; since $A$ already has a canonical $A^{\underline{H}}$-bimodule structure, the trouble lies in defining right and left $\underline{A \# H}$-module structures on $A$. 

Let $\int_H^l$ be the subspace generated by left integrals in $H$. We remind the reader that if $t \in \int_H^l$ then so does $th$ for every $h \in H$. Since 
$\dim \int_H^l = 1$, $th = \alpha(h) t$ for some $\alpha(h) \in k$. This defines a map $\alpha: H \vai k$, which is an algebra morphism.

\begin{lemma} Given  $b \in A$ and $\underline{a\# h} =\sum a(h_{(1)} \cdot \um )\otimes h_{(2)} \in \underline{A \# H}$, the mappings
\[
(\underline{a\# h} ) \seta b = a(h \cdot b) 
\]
and
\[
b \rseta (\underline{a\# h}) =  \sum \alpha(h_{(2)}) S^{-1}(h_{(1)}) \cdot (ba)\]
define left and right $\underline{A \# H}$-module structures  on $A$. Furthermore, if we consider the canonical left and right $A^{\underline{H}}-$ module structures on $A$, then $A$ is both an $A^{\underline{H}} - \underline{A \# H}$ and $\underline{A \# H} - A^{\underline{H}} $ bimodule.
\end{lemma}

\begin{proof}
Let us begin with the left $\underline{A\# H}$ module structure
\[
\begin{array}{rcl}
((\underline{a \# h})(\underline{b \# k}))\seta c & = & 
\sum (\underline{a (h_{(1)} \cdot b)\# h_{(2)}k}) \seta c =\\
& = & (a (h_{(1)} \cdot b)(h_{(2)} k \cdot c) =\\
& = & (a (h \cdot (b(k \cdot c)) =\\
& = & (\underline{a\# h}) \seta (b(k\cdot c))=\\
& = & (\underline{a \# h})\seta ((\underline{b \# k}) \seta c) 
\end{array}
\]
On the right side we have
\begin{eqnarray}
& \, & (a \rseta (\underline{b\# h})) \rseta  (\underline{c \# k})  =\nonumber \\
& = &(\sum \alpha(h_{(2)})(S^{-1}(h_{(1)}) \cdot ab) \rseta (\underline{c \# k}) =\nonumber \\
& = & \sum \alpha(k_{(2)})S^{-1}(k_{(1)}) \cdot ((\sum \alpha (h_{(2)})(S^{-1}(h_{(1)}) \cdot ab) c )) =\nonumber \\
& = & \sum \alpha (h_{(2)}k_{(3)})(S^{-1}(k_{(2)}) \cdot (S^{-1}(h_{(1)}) \cdot ab))(S^{-1}(k_{(1)}) \cdot c) =\nonumber \\
& = & \sum \alpha (h_{(2)}k_{(4)})(S^{-1}(k_{(3)})(S^{-1}(h_{(1)}) \cdot ab))(S^{-1}(k_{(2)}) \cdot \um ) (S^{-1}(k_{(1)}) \cdot c) =\nonumber \\
& = & \sum \alpha (h_{(2)}k_{(3)})(S^{-1}(h_{(1)}k_{(2)}) \cdot ab)(S^{-1}(k_{(1)}) \cdot c) = \nonumber \\
& = & \sum \alpha (h_{(4)}k_{(3)})(S^{-1}(h_{(3)}k_{2}) \cdot ab)(S^{-1}(k_{(1)})S^{-1}(h_{(2)}) h_{(1)} \cdot c) =\nonumber \\
& = & \sum \alpha(h_{(4)}k_{(3)})(S^{-1}(h_{(3)}k_{(2)}) \cdot ab)(S^{-1}(h_{(2)}k_{(1)}) h_{(1)}\cdot c) = \nonumber \\
& = & \sum \alpha(h_{(3)}k_{(2)})(S^{-1}(h_{(2)}k_{(1)}) \cdot (ab ( h_{(1)}\cdot c)) =\nonumber \\
& = & a\rseta (\sum \underline{b(h_{(1)} \cdot c) \# h_{(2)} k} ) =\nonumber \\
& = & a \rseta ((\underline{b\# h}) (\underline{c \# k}) ) .\nonumber
\end{eqnarray}
Hence $A$ is a left and right $\underline{A \# H}$ module. 

 $A$ is an $A^{\underline{H}} - \underline{A \#H} - $ bimodule. Given $a \in A$, $b \in A^{\underline{H}}$ and $\underline{c \# h} \in \underline{A \# H}$, we have
\[
\begin{array}{rcl}
(ba) \rseta (\underline{c \# h}) & = & \sum \alpha(h_{(2)}) S^{-1}(h_{(1)}) \cdot (bac) = \\
& = &  \sum \alpha(h_{(3)}) (S^{-1}(h_{(2)}) \cdot b)(S^{-1}(h_{(1)}) \cdot ac) = \\
& = &  \sum \alpha(h_{(3)}) b (S^{-1}(h_{(2)}) \cdot \um )(S^{-1}(h_{(1)}) \cdot ac) = \\
& = &  \sum \alpha(h_{(2)}) b (S^{-1}(h_{(1)}) \cdot ac) \\
& = & b (a \rseta (\underline{c \# h}))
\end{array}
\]
and 
\[
\begin{array}{rcl}
((\underline{c \# h})\seta a) b & = & c(h \cdot a)b =
 \sum c(h_{(1)} \cdot a) (h_{(2)} \cdot \um) b =\\
& = & c (h \cdot (ab)) = (\underline{c \# h}) \seta (ab).
\end{array}
\]
\end{proof}

For the Morita context we define the maps 
\[
\begin{array}{rcl}
[  \cdot , \cdot ]: A \otimes_{A^{\underline{H}}} A & \vai &  \underline{A \# H}\\
a \otimes b & \mapsto & [a,b] = atb  
\end{array}\]
and 
\[
\begin{array}{rcl}
\langle  \cdot , \cdot  \rangle: A \otimes_{\underline{A \# H}} A & \vai & A^{\underline{H}} \\
a \otimes b & \mapsto &  \langle a,b \rangle = \tr(ab) = t \cdot ab 
\end{array}
\]
Remember that $atb = (\underline{a \# 1_H}) (\underline{\um \# t} ) (\underline{b \# 1_H})$ is indeed an element of $\underline{A \# H}$. 

We must check that these maps are well-defined, i.e., that $[,]$ is $A^{\underline{H}}$-balanced and $\langle, \rangle$ is $\underline{A\# H}$-balanced. Both are clearly $k$-linear maps from $A \otimes_k A$ to $\underline{A \# H}$ and $A^{\underline{H}}$ respectively. Then, we need to check only  whether these maps are balanced or not. 

First for the map $[ \cdot , \cdot ]$: if $a,b \in A$ and $c \in A^{\underline{H}}$ then
\[
\begin{array}{rcl}
[a,cb] & = & atcb = (\underline{a \# 1_H}) (\sum \underline{(t_{(1)} \cdot c) \# t_{(2)}})
(\underline{b \# 1}) =\\
& = & (\underline{a \# 1_H}) (\sum \underline{c (t_{(1)} \cdot \um) \# t_{(2)}})
(\underline{b \# 1_H}) =\\
& = & (\underline{a \# 1_H}) (\underline{c \# 1_H}) (\underline{\um \# t}) 
(\underline{b \# 1_H}) = actb = [ac,b]
\end{array}
\]

Now, for the second map, $\langle  \cdot , \cdot  \rangle$: if $a,b\in A$ and 
$\underline{c\# h}\in \underline{A\# H}$ then
\[
\begin{array}{rcl}
& & \langle a \rseta (\underline{c \# h}), b \rangle = \\
& = & t \cdot (( \sum \alpha(h_{(2)}) S^{-1}(h_{(1)}) \cdot (ac)) b ) =\\
& = & \sum \alpha(h_{(2)}) (t_{(1)} \cdot (S^{-1}(h_{(1)}) \cdot ac))(t_{(2)} \cdot  b ) =\\
& = & \sum \alpha(h_{(2)}) (t_{(1)} S^{-1}(h_{(1)}) \cdot ac)(t_{(2)} \cdot \um )
(t_{(3)} \cdot  b ) =\\
& = & \sum \alpha(h_{(2)}) (t_{(1)} S^{-1}(h_{(1)}) \cdot ac)(t_{(2)} \cdot  b ) \\
& = & \sum \alpha(h_{(3)}) (t_{(1)} S^{-1}(h_{(2)}) \cdot ac)(t_{(2)} \epsilon (h_{(1)}) \cdot  b ) \\
& = & \sum \alpha(h_{(4)}) (t_{(1)} S^{-1}(h_{(3)}) \cdot ac)(t_{(2)}S^{-1}(h_{(2)})h_{(1)} \cdot  b ) =\\
& =& \sum \alpha (h_{(3)})(t S^{-1} (h_{(2)})) \cdot (ac(h_{(1)} \cdot b)) =\\
& =& \sum (t h_{(3)} S^{-1} (h_{(2)})) \cdot (ac(h_{(1)} \cdot b)) =\\
& =& \sum (t \epsilon (h_{(2)})) \cdot (ac(h_{(1)} \cdot b)) =\\
& = & t \cdot (ac (h \cdot b)) = \langle a, (\underline{c \# h}) \seta b \rangle .
\end{array}
\]
Therefore the maps $[\cdot ,\cdot ]$ and $\langle \cdot ,\cdot \rangle$ are well defined.\\

\begin{thm} $(\underline{A \# H}, A^{\underline{H}}, {}_{\underline{A \# H}} A_{A^{\underline{H}}}, {}_{A^{\underline{H}}} \recuo A_{\underline{A \# H}}, [ \cdot , \cdot ], \langle \cdot , \cdot  \rangle)$ is a Morita context.
\end{thm}

\begin{proof}
We must check that both are bimodule maps. From Lemma \ref{trace} it follows that $\langle, \rangle$ is a bimodule map, because $\tr$ is an $A^{\underline{H}}-A^{\underline{H}}-$ bimodule mapping. For the other map, let $a,b \in A$ and $\underline{c \# h} \in \underline{A \# H}$. On the left we have
\[
\begin{array}{rcl}
[(\underline{c \# h}) \seta a, b] & = & c (h \cdot a) t b = 
(\underline{c (h \cdot a) \# t}) (\underline{b \# 1_H}) =\\
& = & \sum (\underline{c (h_{(1)} \cdot a) \# \epsilon (h_{(2)}) t}) (\underline{b \# 1_H}) =\\
& = & \sum (\underline{c (h_{(1)} \cdot a) \# h_{(2)} t}) (\underline{b \# 1_H}) =\\
& = & (\underline{c \# h}) (\underline{a \# t}) (\underline{b \# 1_H}) = (\underline{c \# h}) [a,b].
\end{array}
\]

For the right side, 
\[
\begin{array}{rcl}
& & [a, b \rseta (\underline{c \# h}) ] = \\
& = & [a, \sum \alpha (h_{(2)}) S^{-1} (h_{(1)}) \cdot (bc)]  =\\
& = & (\underline{a \# t}) 
(\sum \underline{\alpha (h_{(2)}) S^{-1} (h_{(1)}) \cdot (bc) \# 1_H})  =\\
& = & \sum (\underline{a \# th_{(2)}}) (\underline{S^{-1} (h_{(1)}) \cdot (bc) \# 1_H})  =\\
& = & (\underline{a \# 1_H}) 
(\sum \underline{t_{(1)}h_{(2)} \cdot (S^{-1} (h_{(1)}) \cdot (bc)) \# t_{(2)}h_{(3)}}) =\\
& = &  (\underline{a \# 1}) 
(\sum \underline{(t_{(1)}h_{(2)} S^{-1} (h_{(1)}) \cdot (bc)) (t_{(2)}h_{(3)} \cdot \um )\# t_{(3)}h_{(4)}})=\\
& = & (\underline{a \# 1}) 
(\sum \underline{(t_{(1)} \cdot (bc)) (t_{(2)}h_{(1)} \cdot \um )\# t_{(3)}h_{(2)}}) =\\
& = &  (\underline{a \# 1}) 
(\sum \underline{(t_{(1)} \cdot ((bc) (h_{(1)} \cdot \um )))\# t_{(2)}h_{(2)}}) =\\
& = & (\underline{a \# 1}) (\underline{\um \# t})(\underline{b \# 1_H})
(\underline{c \# h } )= [a,b](\underline{c \# h }) .
\end{array}
\]

Now we must check the ``associativity'' of the brackets, i.e,  $[a,b] \seta c = 
a \langle b,c\rangle$ and $a \rseta [b,c] = \langle a,b\rangle c$, and we are done. The first is straightforward, and the second is
\[
\begin{array}{rcl}
& & a \rseta [b,c] = \\
& = & a \rseta (\sum \underline{b(t_{(1)} \cdot c) \# t_{(2)}}) =\\
& = & \sum \alpha(t_{(3)})S^{-1}(t_{(2)}) \cdot (ab (t_{(1)} \cdot c)) =\\
& = & \sum \alpha(t_{(5)})(S^{-1}(t_{(4)}) \cdot ab) (S^{-1}(t_{(3)}) \cdot \um ) (S^{-1}(t_{(2)}) t_{(1)} \cdot c) =\\
& = & \sum \alpha(t_{(3)})(S^{-1}(t_{(2)}) \cdot ab) (S^{-1}(t_{(1)}) \cdot \um )  c =\\
& = & \sum \alpha(t_{(2)})(S^{-1}(t_{(1)}) \cdot ab) c =\\
& = & ((\sum \alpha(t_{(2)})(S^{-1}(t_{(1)})) \cdot (ab)) c =\\
& = & (t \cdot (ab)) c =\\
& = & \langle a,b \rangle c ,
\end{array}
\]
where we used $\sum (\alpha(t_{(2)})(S^{-1}(t_{(1)}))  = t$, which follows from $S(t) = \sum \alpha(t_{(2)})t_{(1)}$ (See the reference \cite{radford} for a proof of this result).
\end{proof}

It is worth to mention that in reference \cite{caen06}, the authors considered a partial coaction of a Hopf algebra $H$ on a unital algebra $A$ and established a Morita context between the invariant subalgebra $A^{CoH}$ and the dual smash $\underline{\#}(H,A)\cong {}^*(\underline{A\otimes H})$ but in that case the modules were $A$ as a $A^{\underline{CoH}} -\underline{\#} (H,A)$ bimodule and $Q$ as a $\underline{\#} (H,A) -A^{\underline{CoH}}$ bimodule, where 
\[
Q=\{ q\in {}^* (\underline{A \otimes H}) | c_{(1)} q(c_{(2)})=q(c)\rho (\um ) , \; \forall c \in (\underline{A \otimes H}) \}.
\]

\section{Partial Hopf Galois theory}

In this section, we give some necessary and sufficient conditions for the Morita context $(A^{\underline{H}}, \underline{A\# H}, A, \langle \cdot , \cdot \rangle , [ \cdot , \cdot ] )$ to be strict. This envolves a generalization for the partial case of the concept of a Hopf-Galois extension of algebras. First, given a partial coaction $\prho :A\rightarrow A\otimes H$, we can define an $A-A$ bimodule structure on $A\otimes H$: the left $A$ module structure is given by the multiplication and the right $A$ module structure is given by $(a\otimes h)b=\sum ab^{[0]}\otimes hb^{[1]}$. In what follows, let us consider a sub $A-A$ bimodule of $A\otimes H$ defined by 
\[
\underline{A\otimes H}=(A\otimes H)\um =
\{ \sum  a\um^{[0]} \otimes h\um^{[1]} |\, a\in A, \; h\in H \} .
\]

\begin{defi}
(Partial Hopf Galois Extension, \cite{caen06}) Let $(A, \prho) $ be a partial right $H$-comodule algebra. The extension $A^{\underline{coH}} \subset A$ is partial $H$-Hopf Galois if the canonical map $\can: A \otimes_{A^{\underline{coH}}} A \vai \underline{ A \otimes H }$, given by $\can(a \otimes b) = (a \otimes 1) \prho(b) = \sum a b^{[0]} \otimes b^{[1]}$, is a bijective $A-A$ bimodule morphism.
\end{defi}

\begin{lemma} \label{coH} 
If $H$ is a finite dimensional Hopf algebra, $A$ is a partial $H$-module algebra and $\prho: A \vai A \otimes H^\ast$ is the induced partial $H^\ast$-comodule structure on $A$, then $A^{\underline{H}} = A^{\underline{coH^\ast}}$.
\end{lemma}
\begin{proof}
Let $\{h_i\}_{i=1}^n$ be a basis of $H$ and let $\{h_i^{\ast}\}_{i=1}^n$ be the dual basis. If $a \in A^{\underline{H}}$ then 
\[
\prho(a) = \sum (h_i \cdot a) \otimes h_i^\ast = \sum a(h_i \cdot \um) \otimes h_i^\ast
= a \prho (\um)
\]
(for $\prho(\um) = \sum h_i \cdot \um \otimes h_i^\ast$)

Conversely, if $a \in A^{\underline{coH^\ast}}$ then $\prho(a) = \sum a \um^{[0]} \otimes \um^{[1]}$, and therefore
\[
h \cdot a = (Id\otimes \mbox{ev}_{h})\prho (a) = \sum a \um^{[0]} \um^{[1]}(h) = a (h \cdot \um).
\]
\end{proof}

This last lemma says that when $H$ is finite dimensional, we may consider on a partial $H$-module algebra $A$ the induced structure of partial $H^{\ast}$-comodule algebra and, since $A^{\underline{H}} = A^{\underline{coH^\ast}}$, we get the map \\
$\can: A \otimes_{A^{\underline{H}}} A \vai \underline{A \otimes H^\ast } $.

\begin{thm} \label{kt}
Let $H$ be a finite-dimensional Hopf algebra, $0 \neq t \in \int_H^l$, and let $A$ be a partial $H$-module algebra such that the canonical map $\can: A \otimes_{A^{\underline{H}}} A \vai \underline{A \otimes H^\ast }$ is surjective. Then 
\begin{itemize} 
\item[(i)] There exist $a_1, \ldots, a_n$ and $b_1, \ldots, b_n$ in $A$ such that $\phi_i: A \vai A^{\underline{H}}$ given by $\phi_i(a) = t \cdot b_ia$ is an $A^{\underline{H}}$ module map, and  $a = \sum a_i \phi_i(a)$ for each $a \in A$; hence $\{a_i\}_{i=1}^n$ is a projective basis over $A^{\underline{H}}$ and $A$ is a finitely generated projective $A^{\underline{H}}$ right module.
\item[(ii)] $\can$ is bijective.
\end{itemize}
\end{thm}

\begin{proof}
(i) Consider the canonical isomorphism \cite{susan} induced by the nonzero left integral $t\in \int_H^l$ given by
\[
\begin{array}{rcl} \theta :H^\ast & \rightarrow & H\\
                            f     & \mapsto & \theta (f)=t\hitby f
\end{array}
\]
where $t\hitby f=\sum f(t_{(1)})t_{(2)}$. Then, as $\theta$ is surjective, there exists $T\in H^\ast$ such that $1_H =t\hitby T$. Using also the surjectivity of the canonical map $\beta$, there exist $a_1$, $\ldots$, $a_n$ and $b_1$, $\ldots$, $b_n$, elements of $A$ such that
\[
\um^{[0]}\otimes T\um^{[1]} =\can (\sum_{i=1}^n a_i \otimes_{A^{\underline{H}}} b_i ).
\]
Now, consider any element $a\in A$, then
\begin{eqnarray}
a=1_H\cdot a &=& (t\hitby T)\cdot a =(t\hitby T) (\um a)=\nonumber\\
&=& \sum \um^{[0]}a^{[0]} (\um^{[1]} a^{[1]} (t\hitby T)) =\nonumber\\
&=& \sum \um^{[0]} a^{[0]} T(t_{(1)})(\um^{[1]}a^{[1]} (t_{(2)}) )=\nonumber\\
&=& \sum \um^{[0]} a^{[0]} T(t_{(1)}) \um^{[1]}(t_{(2)})a^{[1]}(t_{(3)}) =\nonumber\\
&=& \sum \um^{[0]} a^{[0]} (T\um^{[1]}(t_{(1)}))a^{[1]}(t_{(2)}) =\nonumber\\
&=& \sum_{i=1}^n \sum a_i b_i^{[0]} a^{[0]} b_i^{[1]}(t_{(1)})a^{[1]}(t_{(2)})=\nonumber\\
&=& \sum_{i=1}^n \sum a_i b_i^{[0]} a^{[0]} (b_i^{[1]}a^{[1]}(t))=\nonumber\\
&=& \sum_{i=1}^n \sum a_i (t\cdot b_i a).\nonumber
\end{eqnarray}

(ii) the proof of this item follows the same steps as in \cite{susan} Theorem 8.3.1.
\end{proof}

\begin{prop}
Let $H$ be a finite dimensional Hopf algebra and $A$ a partial $H$-module algebra. Then $A^{\underline{H}}\cong \mbox{End}({}_{\underline{A\# H}} A)^{op}$ as algebras.
\end{prop}

\begin{proof}
Let $\sigma: A^{\underline{H}}\vai \mbox{End}({}_{\underline{A\# H}} A)^{op}$ be the map that takes $a$ to the endomorphism $\sigma(a): b \mapsto ba$. Then $\sigma$ is a algebra map and, if $f \in \mbox{End}({}_{\underline{A\# H}} A)$ and $a \in A$, then 
$$f(a) = f(a \um) = f ((\underline{a \# 1}) \cdot \um) = (\underline{a\# 1}) \cdot 
f(\um ) = a f(\um).$$

Hence $f = \sigma_{f(\um)}$;  besides, $f(\um )$ is indeed in $A^{\underline{H}}$, since
$$h \cdot f(\um ) = \sum (\underline{(h_{(1)} \cdot \um )  \# h_{(2)} }) \cdot f(\um ) = f((h \cdot \um ) \um ) = (h \cdot \um ) f(\um). $$
\end{proof}  

Finally, the next result provides a relation between the surjectivity of the canonical map, the extension being Hopf-Galois, and the surjectivity of the map $[,]$ of the Morita context.

\begin{thm} 
Let $H$ be a finite dimensional Hopf algebra with a nonzero integral $t$, and let $A$ be a (left) partial $H$-module algebra. Then, the following affirmations are equivalent:
\begin{itemize}
\item[(i)]  The canonical map $\beta :A\otimes_{A^{\underline{H}}} A \rightarrow \underline{A\otimes H^* }$, given by $\beta (\sum x_i \otimes y_i ) =\sum x_i y_i^{[0]} \otimes y_i^{[1]}$ is surjective.
\item[(ii)] The algebra $A$ is a finitely generated projective right $A^{\underline{H}}$ module and $A^{\underline{H}} \subset A$ is a partial $H^\ast$-Galois extension.
\item[(iii)] The algebra $A$ is a finitely generated projective right $A^{\underline{H}}$ module, and the map $[\cdot , \cdot ]:A\otimes_{A^{\underline{H}}} A \rightarrow \underline{A\# H}$ is surjective
\end{itemize}
\end{thm}

\begin{proof}
(i) $\Rightarrow $ (ii) This is basically the content of the Theorem \ref{kt}

(ii) $\Rightarrow$ (iii)

Let $\theta: H^\ast \vai H$ be the $H$-module isomorphism $\varphi \mapsto t \hitby \varphi = \sum \varphi(t_1)t_2$, and let
$\can$ be the canonical map ($\can(a \otimes b) = \sum ab^{[0]} \otimes b^{[1]}$). Then $[a,b] = (I \otimes \theta) \can (a \otimes b)$. In fact, 
\[
\begin{array}{rcl}
(I \otimes \theta) \can (a \otimes b) & = & ab^{[0]} \otimes \theta (b^{[1]}) \\
& = & \sum ab^{[0]} \otimes  t \hitby (b^{[1]}) \\
& = & \sum  ab^{[0]} \otimes  (b^{[1]}) (t_1) t_2 \\
& = & \sum a (t_1 \cdot b) \otimes t_2 \\
& = & \sum (a (t_1 \cdot \um) \otimes t_2 )(b\otimes 1_H) \\
& = & (\underline{a \# t})(\underline{b \# 1_H})\\
& = & (\underline{a \# 1_H })(\underline{\um \# t})(\underline{b \# 1_H})\\
& = & [a,b]
\end{array}
\]
This formula assures that $[\cdot , \cdot ]$ is surjective if, and only if $\beta$ is surjective. Hence, if (ii) holds then $\can$ is bijective, in particular, $\beta$, therefore $[,]$ is onto.

(iii)$\Rightarrow$ (i) Note that, the formula $[\cdot ,\cdot ]=(I\otimes \theta)\beta$ implies that the surjectivity of the map $[\cdot , \cdot ]$ assures the surjectivity of $\beta$.
\end{proof}

Note that the previous theorem gives us a necessary and sufficient condition for the Morita context $(A^{\underline{H}}, \underline{A\# H}, A, \langle \cdot , \cdot \rangle , [ \cdot , \cdot ] )$ to be strict. In fact, assuming that $H$ is a finite dimensional Hopf algebra with a nonzero integral $t$ and that the map $\tr :A\rightarrow A^{\underline{H}}$, defined as $\tr (a)=t\cdot a$, is surjective, we have that the Morita context is strict if, and only if, the extension $A^{\underline{H}} \subset A$ is partial $H^{\ast}$-Hopf Galois.

\section*{Aknowledgements}

This work is dedicated to Prof. Miguel Ferrero, whose contribution to the development of algebra in Brazil is remarkable.

The authors would like to thank to Piotr Hajac and Joost Vercruysse for fruitfull discussions and suggestions.

\end{document}